\documentclass[1p]{elsarticle}
\usepackage{amssymb}
\usepackage{amsfonts}

\usepackage[polish, english]{babel}
\usepackage{polski}

\newtheorem{theorem}{Theorem}
\newtheorem{conjecture}[theorem]{Conjecture}
\newtheorem{corollary}[theorem]{Corollary}
\newtheorem{lemma}[theorem]{Lemma}

\newtheorem{claim}{Claim}

\newproof{pf}{Proof}

\begin{document}

\title{On conflict-free proper colourings of graphs without small degree vertices}

\author[agh]{Mateusz Kamyczura}
\ead{kamyczuram@gmail.com}

\author[agh]{Jakub Przyby{\l}o}
\ead{jakubprz@agh.edu.pl}

\address[agh]{AGH University of Science and Technology, Faculty of Applied Mathematics, al. A. Mickiewicza 30, 30-059 Krakow, Poland}

\begin{abstract}
A proper vertex colouring of a graph $G$ is referred to as conflict-free if in the neighbourhood of every vertex some colour appears exactly once, while it is called $h$-conflict-free if there are at least $h$ such colours for each vertex of $G$. The least numbers of colours in such colourings of $G$ are denoted $\chi_{\rm pcf}(G)$ and $\chi_{\rm pcf}^h(G)$, respectively. It is known that  $\chi_{\rm pcf}^h(G)$ can be as large as $(h+1)(\Delta+1)\approx \Delta^2$ for graphs with maximum degree $\Delta$ and $h$ very close to $\Delta$.  We provide several new upper bounds for these parameters for graphs with minimum degrees $\delta$ large enough and $h$ detached from $\delta$. In particular we show that $\chi_{\rm pcf}^h(G)\leq (1+o(1))\Delta$ if $\delta\gg\ln\Delta$ and $h\ll \delta$, and that $\chi_{\rm pcf}(G)\leq \Delta+O(\ln \Delta)$ for regular graphs. These specifically refer to the conjecture of Caro, Petru\v{s}evski and \v{S}krekovski that $\chi_{\rm pcf}(G)\leq \Delta+1$ for every connected graph $G$ of maximum degree $\Delta\geq 3$, towards which they proved that $\chi_{\rm pcf}(G)\leq \left\lfloor\frac{5\Delta}{2}\right\rfloor$ if $\Delta\geq 1$.
\end{abstract}

\begin{keyword}
conflict-free colouring \sep conflict-free chromatic number \sep $h$-conflict-free chromatic number \sep odd colouring
\end{keyword}

\maketitle

\section{Introduction and notation}
A conflict-free colouring is a vertex colouring of a graph such that for each vertex there is a colour appearing exactly once in its neighbourhood. This type of colouring was first defined in the environment of hypergraphs, where it was required that each hyperedge contains a unique colour, i.e. a colour appearing exactly once. Such concept was originally motivated by channel assignment in wireless networks, where colours represent available frequencies, which potentially interfere, and thus one should always have in range a transmitter with a unique frequency associated, see e.g.~\cite{EvenEtAl,SmorodinskyPhd,SmorodinskyApplications}. In 2009 Pach and Tardos~\cite{PachTardos} provided a general upper bound of order close to optimal for the related graph invariant. Its precise order of magnitude was later settled down by Bhyravarapu, Kalyanasundaram and Mathew~\cite{Hindusi} and by Glebov, Szab\'o and Tardos~\cite{GlebovEtAl}. 
See also~\cite{Hindusi2,KellerEtAl,KostochkaEtAl} for other related results.

A variant of the parameter above for proper colourings was later introduced by Fabrici, Lu\v{z}ar, Rindo\v{s}ov\'a and Sot\'ak, who provided a series of results concerning planar graphs, see~\cite{Fab} for details. Consider a proper vertex colouring of a graph $G=(V,E)$ without isolated vertices.
We shall say that a vertex $v\in V$ has a unique neighbour if in the neighbourhood of $v$ there exists a vertex $u$ whose colour is unique among all colours of the neighbours of $v$. 
A proper vertex colouring of $G$ with $k$ colours such that each $v\in V$ has a unique neighbour 
is referred to as a conflict-free $k$-colouring (or a conflict-free $[k]$-colouring if $[k]:=\{1,2,\ldots,k\}$ is the set of used colours). The least $k$ for which $G$ is conflict-free $k$-colourable is called the conflict-free chromatic number of $G$. We denote it by $\chi_{\rm pcf}(G)$. In~\cite{Fab} it was specifically  proven that $\chi_{\rm pcf}(G)\leq 8$ for planar graphs and conjectured that $\chi_{\rm pcf}(G)\leq 6$ in such a case. The latter was confirmed in~\cite{CPS} for planar graphs of girth at least $7$  by Caro, Petru\v{s}evski and \v{S}krekovski, who also obtained several other related results. 
 Notably, they proved the following general upper bound.
\begin{theorem}[\cite{CPS}]
For every connected graph $G$ of maximum degree $\Delta\geq 1$,
$$\chi_{\rm pcf}(G)\leq \left\lfloor\frac{5\Delta}{2}\right\rfloor.$$
\end{theorem}
They also boldly conjectured that if only the maximum degree of $G$ is not very small, i.e. at least $3$, then something much stronger holds.
\begin{conjecture}[\cite{CPS}]\label{ConjectureCPS}
For every connected graph $G$ of maximum degree $\Delta\geq 3$,
$$\chi_{\rm pcf}(G)\leq \Delta+1.$$
\end{conjecture}
See e.g.~\cite{EUN,EUN2,DWC-CHL,RH,CHL} for further and related results. In~\cite{DWC-CHL} Cranston and Liu in particular provided  an upper bound $\chi_{\rm pcf}(G)\leq 1.6550826\Delta+o(\Delta)$, which holds also in a more general list setting.
Moreover, in~\cite{EUN2},  a natural extension of the concept of conflict-free colourings, where each vertex is required to have at least $h$ unique colours,  rather than just one, in its neighbourhood, was considered. Such proper vertex colouring of a graph $G$ with $k$ colours is called an $h$-conflict-free $k$-colouring of $G$. The least $k$ for which $G$ is $h$-conflict-free $k$-colourable is in turn called the $h$-conflict-free chromatic number of $G$ and we denote it by $\chi_{\rm pcf}^h(G)$ (hence $\chi_{\rm pcf}^1(G) = \chi_{\rm pcf}(G)$). Cho, Choi, Kwon and Park~\cite{EUN2} specifically observed 
that for infinitely many $h$, there exists a graph $G$ for which $\chi_{\rm pcf}^h(G) = (h+1)(\Delta-1)$. Such graphs may be constructed with aid of orthogonal Latin squares. In these examples however, $h$ is very close to $\Delta$, which implies that $\chi_{\rm pcf}^h(G) \approx \Delta^2$ for them. 
Moreover, the colourings imposed by these examples must in fact be square colourings (called also $2$-distance colourings), which attribute distinct colours to vertices at distance at most $2$ -- the least number of colours in these define the well studied graph invariant, the $2$-chromatic number of a graph, $\chi_2(G)$ (which equals the chromatic number of the square of the given graph $G$), for which infinitely many examples of graphs with $\chi_2(G) \geq \Delta^2-O(\Delta)$ are known, see e.g.~\cite{BDMP,DistChrSurvey}. Note that $\chi_{\rm pcf}^h(G)$ may actually be regarded as a relaxation of $\chi_2(G)$, while the two parameters coincide for $\Delta$-regular graphs and $h=\Delta$ (as well as $h=\Delta-1$).

In the present note we show that Conjecture~\ref{ConjectureCPS} holds asymptotically if the minimum degree of $G$ is large enough, i.e. $\delta\gg\log\Delta$. In fact we prove that the same conclusion holds also for $h>1$, if only $h\ll \delta$, cf. Corollary~\ref{1+o(1)Corollary} below.
For regular graphs we moreover specify that $\chi_{\rm pcf}(G)\leq \Delta+O(\log\Delta)$, and similarly $\chi_{\rm pcf}^h(G)\leq \Delta+O(\log\Delta)$ even for $h=O(\log\Delta)$, cf. Corollary~\ref{RegularCFCorollary}.
We also provide several other observation based on our main probabilistic argument. 

Let us finally remark that the conflict-free colouring may also be regarded as a direct strengthening of the recently intensively studied odd colouring, which is a proper vertex colourings such that for each vertex $v$ there is a colour appearing an odd number of times in the neighbourhood of $v$. The least number of colours admitting such a colouring of a graph $G$ is denoted $\chi_{\rm odd}(G)$. Our main outcome, Theorem~\ref{MainComplicated}, thus provides strengthenings of some known results  concerning odd colourings, in particular in reference to the general upper bound $\chi_{\rm odd}(G)\leq 2\Delta$ for  every connected graph $G$ of maximum degree $\Delta$ except $C_5$, attained by Caro,  Petru\v{s}evski and \v{S}krekovs~\cite{ODD3}, who also conjectured that in fact $\chi_{\rm odd}(G)\leq \Delta+1$ for every connected graph $G$ of maximum degree $\Delta\geq 3$, cf. Conjecture~\ref{ConjectureCPS} above. See also e.g.~\cite{ODD3,ODD5,EUN2,ODD2,ODD6,RH,ODD4,ODD1} for other results concerning this graph invariant.

In the next section we briefly recall probabilistic tools applicable within our arguments.
In the third section we present our main result, whose formulation might seem slightly  inapparent and technical. 
The last section includes several more self-evident conclusions and a few comments, concerning in particular directions for possible further research.

\section{Basic probabilistic tools}
We shall use the basic symmetric variant of the Lov\'{a}sz Local Lemma, which can be found e.g. in \cite{LLL}. 

\begin{lemma}[Lov\'asz Local Lemma]\label{LLL}
Let $\Omega$ be a finite family of events in any probability space. Suppose that every event $A \in \Omega $ is mutually independent of a set of all the other events in $\Omega$ but at most $D$, and that $\mathbf{Pr}(A) \leq p$ for each $A \in \Omega$. If 
$$ep(D + 1) \leq 1,$$
then $\mathbf{Pr}(\bigcap_{A \in \Omega} \overline{A}) > 0$. 
\end{lemma}

We shall also need the Chernoff Bounds, see e.g.~ \cite{RandOm}. 
\begin{lemma}[Chernoff Bounds]\label{Ch} Let $X= \sum_{i=1}^n X_i$, where $X_i=1$ with probability $p_i$, $X_i = 0$ with
probability $1 - p_i$ and all $X_i$ are independent. Then:
\begin{itemize}
    \item $ \mathbf{Pr}(X < (1-\delta)\mathbf{E}(X))< \exp(-\delta^2\mathbf{E}(X)/2) ~~ for~all  ~~ 0<\delta \leq 1,$
    \item $  \mathbf{Pr}(X > (1+\delta)\mathbf{E}(X))< \exp(-\delta^2\mathbf{E}(X)/3) ~~  for~all  ~~ 0<\delta \leq 1.$
\end{itemize}
\end{lemma}
It is straightforward to notice that the Chernoff Bounds in Lemma~\ref{Ch} may also be applied for appropriate estimates on $\mathbf{E}(X)$, the first of these if $\mathbf{E}(X)\geq a>0$ and the second for $\mathbf{E}(X)\leq a$, with $\mathbf{E}(X)$ replaced by $a$ in the both inequalities.

\section{Main result}
\begin{theorem}\label{MainComplicated}
There exists $\Delta_0$ such that  
$\chi^h_{\rm pcf}(G)< \Delta(1+\max\{30\frac{h}{\delta},600\frac{\ln\Delta}{\delta}\})$
for every graph $G$ with maximum degree $\Delta\geq \Delta_0$,
minimum degree $\delta \geq 1500 \ln\Delta$ and $h\leq \frac{\delta}{75}$.
\end{theorem}

\begin{pf}
Suppose $G=(V,E)$ and $h\geq 1$ comply with the theorem's assumptions.
Note we do not specify $\Delta_0$. Whenever needed we shall however assume that $\Delta$ is large enough so that all explicit inequalities below hold. Let us fix any linear ordering on $V$. Note it naturalny induces the ordering in any subset of vertices of $G$. Set

\begin{equation}\label{h'_and_m_definitions}
h' : = \max\left\{h,20\ln\Delta\right\}~~~~~~{\rm and} ~~~~~~ m : =  \left\lceil 30\frac{\Delta}{\delta}h' \right\rceil.
\end{equation}
Note that 
\begin{equation}\label{h'ineq}
h'\leq \frac{\delta}{75}.
\end{equation}
Suppose that initially all vertices of $G$ are uncoloured. 
We shall first colour some of them by repeating $m$ times a similar random procedure specified below.

For technical reasons, in each of $m$ rounds of the procedure, we shall however be applying it to a certain supergraph of $G$ (though eventually we shall ignore all additional vertices and their colours). Let $G_1=G$. For $i=2,3,\ldots,m$, we construct $G_i$ from $G_{i-1}$ as follows. For each consecutive vertex $v\in V$ (where by $d(v)$ we shall always mean $d_G(v)$) add $\lceil d(v)/2\rceil$ new (uncoloured) vertices, join them by edges with $v$ and extend our ordering by placing the new vertices at its end (as its largest elements). Denote by $V_i$ the set of vertices of $G_i$. 

For $i=1,2,\ldots,m$ we subsequently proceed as follows:
\begin{enumerate}
    \item[(a)] independently choose every vertex in $V_i$ at random with probability $\frac{1}{\Delta}$; denote the set of all chosen vertices by $C^{\rm ch}_i$; 
    \item[(b)] denote by $C^{\rm iso}_i$ the set of all vertices in $C^{\rm ch}_i$ which have no neighbours in $C^{\rm ch}_i\cap V$;
    \item[(c)] denote all thus far uncoloured vertices in  $C^{\rm iso}_i$ by $C_i$; colour them with $i$.
\end{enumerate}
We shall show that with positive probability each vertex shall have at least $h'\geq h$ uniquely coloured neighbours in the obtained partial $[m]$-colouring of $G$, which is proper due to (b). 

We note that at first glance this procedure may seem slightly wasteful, as we in particular do not admit colouring with $i$ a vertex $v\in C^{\rm ch}_i\cap V$ if it has a neighbour $u\in C^{\rm ch}_i\cap V$, even if $u$ was coloured in some earlier step with colour $j<i$ (hence $u$ cannot be coloured $i$, and thus cannot violate properness of the colouring in round $i$).
Such setup however does not worsen the obtained result, while allowing us to dispose of some problematic dependencies, which would significantly constrict the analysis of our random process.  

For each $v\in V$ and $i\in [m]$ denote the following events:
\begin{itemize}
    \item $A_v^i$: among the first $\lceil d(v)/2\rceil$ neighbours of $v$ in $G_i$ (according to our fixed vertex ordering) which were not coloured in the first $i-1$ rounds of the procedure there is 
    $u$ which was coloured in round $i$ (i.e. $u\in C_i$, which is equivalent to $u\in C^{\rm iso}_i$) and no other neighbour of $v$ in $G$ was chosen in round $i$ (i.e. $w\notin C^{\rm ch}_i$ for each $w\in N_G(v)\smallsetminus\{u\}$);
    \item $A_v$: at least $h'$ events in $\{A_v^1,...,A_v^m\}$ hold;
    \item $B_v$: at most $\frac{d(v)}{2}$ neighbours of $v$ in $G$ were coloured within all $m$ rounds.  
\end{itemize} 
We shall use the Lov\'asz Local Lemma in order to show that $\mathbf{Pr}(\bigcap_{v\in V} (A_v\wedge B_v))>0$, i.e. that there exists a partial proper colouring of $G_m$ such that $A_v$ and $B_v$ are fulfilled for all $v\in V$. Note that such a colouring narrowed down to $G$ would be a partial proper colouring of $G$ within which every vertex $v\in V$ has at least $h'\geq h$ uniquely coloured neighbours -- this follows from the fact that $B_v$ guarantees that in each case the vertex  $u$ from the definition of $A_v^i$ must belong to $V$.
Clearly each $A_v$ and similarly each $B_v$ is determined by random choices committed for vertices at distance at most $2$ from $v$, hence each of $A_v$ and $B_v$ is mutually independent from all other events $A_u$ and $B_u$ except possibly these for which $u$ is at distance at most $4$ from $v$, i.e. all but at most $2\Delta^4$ such events. By Lemma~\ref{LLL} it is thus sufficient to prove the following claim.
\begin{claim}\label{ClaimAvBv}
For each $v\in V$, $\mathbf{Pr}(\overline{A}_v) \leq \Delta^{-5}$ and $\mathbf{Pr}(\overline{B}_v) \leq \Delta^{-5} $.
\end{claim}

\begin{pf} (of Claim \ref{ClaimAvBv}) 
Consider an auxiliary $0$-$1$ random matrix $S=[s_{ij}]_{m\times (d(v)+\Delta)}$ whose every entry is independently set as $1$ with probability $\frac{1}{\Delta}$ or as $0$ with probability $1-\frac{1}{\Delta}$.
Let $X$ denote the number of rows in $S$ each of which contains exactly one $1$ and this $1$ belongs to one of $\lceil \frac{d(v)}{2}\rceil$ initial columns of $S$.
We claim that the probability of $\overline{A}_v$ is bounded above by the probability that $X<h'$. Indeed, for every $i=1,2,\ldots,m$ let us couple the choices for $s_{i,1},\ldots,s_{i,\lceil d(v)/2\rceil}$  with the choices for the first $\lceil\frac{d(v)}{2}\rceil$ uncoloured neighbours of $v$ in $G_i$ within the $i$-th round of our procedure (note that due to definition of $G_i$ there are always sufficiently many such neighbours). Next let us couple choices for possible yet uncoupled neighbours of $v$ in $G$ with the consecutive (at most $\lfloor\frac{d(v)}{2}\rfloor$) entries in the $i$-th row of $S$, and finally couple the entries $s_{i,d(v)+1}, s_{i,d(v)+2},\ldots$ with the choices for yet uncoupled neighbours in $V$ of the first uncoloured (before $i$-th round) neighbour of $v$ (in $G_i$) which was chosen in the $i$-th round of the procedure (if there are any). Note that not all entries in the $i$-th row of $S$ must have been coupled. Nevertheless, the probability of $A_v^i$ is surely not smaller than the probability that the $i$-th row of $S$ contains exactly one $1$ and this $1$ belongs to one of $\lceil \frac{d(v)}{2}\rceil$ initial columns of $S$. Note that the latter equals
 \begin{eqnarray}
 \left\lceil\frac{d(v)}{2}\right\rceil\cdot \frac{1}{\Delta} \cdot \left(1-\frac{1}{\Delta}\right)^{d(v)+\Delta-1}
 &>& \left\lceil\frac{d(v)}{2}\right\rceil\cdot \frac{1}{\Delta} \cdot e^{-\frac{1}{\Delta-1}(d(v)+\Delta-1)}\nonumber\\
&\geq& \left\lceil\frac{d(v)}{2}\right\rceil\cdot \frac{1}{\Delta} \cdot e^{-(2+\frac{1}{\Delta-1})}\nonumber\\
  &>& \frac{1}{15} \cdot \frac{d(v)}{\Delta}\nonumber
\end{eqnarray}
for $\Delta$ large enough.
Hence $\mathbf{E}(X)>\frac{md(v)}{15\Delta}\geq 2h'$, and since all choices in $S$ are independent, by the Chernoff Bound, 
\begin{equation}\label{Av_bound}
\mathbf{Pr}(\overline{A}_v)\leq \mathbf{Pr}(X<h')\leq e^{-\frac{h'}{4}}\leq e^{-5\ln\Delta}=\Delta^{-5}.
\end{equation}

In order to prove that the probability of $\overline{B}_v$ is small it is on the other hand sufficient to prove that the probability that a potentially larger random variable than the one related with $\overline{B}_v$, namely expressing the number of times a neighbour of $v$ in $G$ was chosen at all, exceeds $d(v)/2$ is small. Thus set 
$$Y:=\sum_{i=1}^m|C_i^{\rm ch}\cap N_G(v)|.$$
Note that by~(\ref{h'_and_m_definitions}) and~(\ref{h'ineq}),
$$\mathbf{E}(Y)=\frac{md(v)}{\Delta} < 31\frac{d(v)}{\delta}h' < 0.42d(v).$$ 
Since all choices are independent, by the Chernoff Bound and the fact that $\delta\geq 1500\ln\Delta$ we thus have that
\begin{equation}\label{Bv_bound}
\mathbf{Pr}(\overline{B}_v)\leq \mathbf{Pr}\left(Y>\frac{d(v)}{2}\right)\leq e^{-\frac{-\left(\frac{0.08}{0.42}\right)^2\cdot 0.42d(v)}{3}} < e^{-0.005d(v)}\leq e^{-7.5\ln\Delta}<\Delta^{-5}.
\end{equation}
The claim follows by~(\ref{Av_bound}) and~(\ref{Bv_bound}). \qed
\end{pf}

By Claim~\ref{ClaimAvBv} and the Lov\'asz Local Lemma there exists a partial proper $[m]$-colouring $c$ of $G$ such that every vertex $v\in V$ has at least $h'\geq h$ uniquely coloured neighbours. Let $G'$ be a graph induced by uncoloured vertices in $G$, hence 
$$\Delta(G')\leq \Delta-h' \leq \Delta - 20\ln\Delta < \Delta - 2$$
for $\Delta$ large enough. Clearly there thus exists a proper vertex colouring of $G'$ with less than $\Delta-1$ new colours (not belonging to $[m]$). This combined with $c$ yields a desired $h$-conflict-free colouring of $G$ with less than  $\Delta(1+\max\{30\frac{h}{\delta},600\frac{\ln\Delta}{\delta}\})$ colours.
\qed
\end{pf}

\section{Conclusions and final remarks}
Note we did not put much effort in optimizing multiplicative constants in Theorem~\ref{MainComplicated}.
Hence these could still be somewhat improved in the corollaries below. 
Observe first that since the assumptions in Theorem~\ref{MainComplicated} imply that $30\frac{h}{\delta}\leq 0.4$ and $600\frac{\ln\Delta}{\delta}\leq 0.4$, we obtain the following.
\begin{corollary}\label{1.4Corollary}
There exists $\Delta_0$ such that $\chi^h_{\rm pcf}(G)< 1.4 \Delta$ for every graph $G$ with maximum degree $\Delta\geq \Delta_0$, minimum degree $\delta \geq 1500 \ln\Delta$ and $h\leq \frac{\delta}{75}$.
\end{corollary}
The upper bound for $\chi^h_{\rm pcf}(G)$ which follows from Theorem~\ref{MainComplicated} may also be improved by strengthening multiplicative constants in assumptions on $\delta$ and $h$.
\begin{corollary}
For every $\varepsilon>0$ there exist $C,\gamma$ such that $\chi^h_{\rm pcf}(G)< (1+\varepsilon) \Delta$ for every graph $G$ with maximum degree $\Delta$, minimum degree $\delta \geq C \ln\Delta$ and $h\leq \gamma\delta$.
\end{corollary}
Via strengthening the mentioned assumptions by any relevant factor, changing their order, 
one may moreover derive from Theorem~\ref{MainComplicated} the following, possibly most important corollary, which implies that Conjecture~\ref{ConjectureCPS} of Caro, Petru\v{s}evski and \v{S}krekovski 
holds asymptotically if only the order of $\delta$ exceeds $\log\Delta$. 
\begin{corollary}\label{1+o(1)Corollary}
Let $G$ be a graph with minimum and maximum degrees $\delta,\Delta$, respectively. Then
$\chi^h_{\rm pcf}(G) < (1+o(1)) \Delta$ if $\delta \gg \ln\Delta$ and $h\ll \delta$.
\end{corollary}
Further, in the case of regular graphs specifically and relatively small $h$ one may also derive an upper bound which is apparently close to $\Delta+1$ postulated by Conjecture~\ref{ConjectureCPS}.
\begin{corollary}\label{RegularCFCorollary}
There exists $\Delta_0$ such that $\chi^h_{\rm pcf}(G)< \Delta + 600\ln \Delta $ for every $\Delta$-regular graph $G$ with $\Delta\geq \Delta_0$ and  $h\leq 20 \ln \Delta$.
\end{corollary}
As mentioned in the introduction, all bounds above transfer also directly to environment of odd colourings, in particular Corollary~\ref{1+o(1)Corollary} confirms asymptotically the correspondent of Conjecture~\ref{ConjectureCPS} concerning $\chi_{\rm odd}(G)$ of Caro, Petru\v{s}evski and \v{S}krekovski for 
$\delta\gg\log\Delta$, while Corollary~\ref{RegularCFCorollary} implies e.g. the following conclusion referring to this conjecture and results from~\cite{ODD3}.
\begin{corollary}\label{RegularODDCorollary}
There exists $\Delta_0$ such that $\chi_{\rm odd}(G)< \Delta + 600\ln \Delta $ for every $\Delta$-regular graph $G$ with $\Delta\geq \Delta_0$.
\end{corollary}

Observe finally that e.g. Corollary~\ref{1.4Corollary} implies an upper bound of order $\Delta$ for $\chi^h_{\rm pcf}(G)$ in the case of $\Delta$-regular graphs and $h\leq\frac{\Delta}{75}$, while on the other hand, by research concerning square colourings, we know that $\chi^h_{\rm pcf}(G)$ can be of order $\Delta^2$ for $h$ very close to $\Delta$. It thus seems very intriguing to investigate the behaviour of extreme values of $\chi^h_{\rm pcf}(G)$ with reference to $h$, and in particular finding a threshold for $h$ in terms of $\Delta$ above which $O(\Delta)$ of colours might not suffice anymore, in particular in the case of regular graphs. 

Clearly, it would also be interesting to get closer to proving Conjecture~\ref{ConjectureCPS}, e.g. via designing an approach providing an upper bound of the form $\Delta+const.$, even in the case restricted to regular graphs exclusively.

\end{document}